\documentclass{amsart}

\usepackage{amssymb}
\usepackage{amsmath}
\usepackage{enumerate}
\usepackage{array}
\usepackage{todonotes}
\usepackage{color}
\usepackage{dsfont}
\usepackage{comment} 
\usepackage[utf8]{inputenc}
\usepackage{hyperref}
\usepackage{blindtext}

\newtheorem{Theorem}{Theorem}
\newtheorem{Lemma}{Lemma}

\newtheorem{Corollary}{Corollary}
\newtheorem{Proposition}{Proposition}

\theoremstyle{remark}
\newtheorem{Remark}{Remark}

\numberwithin{equation}{section}

\newcommand{\q}{\mathfrak{q}}

\newcommand{\mfp}{\mathfrak{p}}
\newcommand{\mfq}{\mathfrak{q}}

\newcommand{\T}{\mathbb{T}}

\newcommand{\Z}{\mathbb{Z}}

\newcommand{\mP}{\mathbb{P}}
\newcommand{\mE}{\mathbb{E}}

\newcommand{\var}{\textrm{var}}

\newcommand{\tilderho}{\widetilde{\rho}_x}

\title[Mixing of...]{Mixing of a generic simple symmetric random walk on the circle}
\author{Klaudiusz Czudek}
\address{Klaudiusz Czudek, Institute of Science and Technology Austria (ISTA), Am Campus 1, 3400 Klosterneuburg}
\email{klaudiusz.czudek@gmail.com}

\subjclass[2020]{ Primary 37A25, 60K37.}

\keywords{mixing, Diophantine approximation, quasiperiodic environment}

\begin{document}

\begin{abstract}
Fix an irrational number $\alpha$ and a real function $\mathfrak{p}$ on the circle with $0<\mathfrak{p}<1$. If a particle is placed at a point $x\in \mathbb R/\mathbb Z$, then in the next step it jumps to $x+\alpha$ with probability $\mathfrak{p}(x)$ and to $x-\alpha$ with probability $1-\mathfrak{p}(x)$. Sinai and Kaloshin proved that if $\mathfrak{p}$ is smooth then the random walk is uniquely ergodic and mixing, unless $\alpha$ is Liouville and $\mathfrak{p}$ is symmetric. 
Unique ergodicity in the general case has been obtained by Conze and Guivarc'h. Here we give an alternative proof of the latter as well as some generic result about mixing, which partially solves a recent open problem.
\end{abstract}

\maketitle

\section{Introduction}
\subsection{Definition of the process}

Fix an irrational number $\alpha\in \mathbb R$, a continuous function $\mathfrak{p} : \mathbb{T} \rightarrow (0,1)$, and consider a Markov process $(X_n)$ with transition kernel
\begin{equation}\label{E:1.1}
p(x, \cdot ) = \mathfrak{p}(x) \delta_{x+\alpha} + \q(x) \delta_{x-\alpha}, \quad p : \T \times \mathcal B (\T ) \to [0,1],
\end{equation}
where  $\mathcal B (\mathbb \T )$ stands for the $\sigma$-algebra of Borel subsets of $\mathbb \T$ and $\q(x)=1-\mathfrak{p}(x)$, $x\in \T$. For any Borel probability measure $\mu$ on the circle, the symbol $\mP_\mu$ denotes the conditional probability given the process (\ref{E:1.1}) is launched with the initial distribution $\mu$. A Borel probability measure $\mu$ is called \textit{stationary} if for every   $\varphi$ from the space of continuous functions $C(\T)$ on $\T$ we have $\mE_\mu\varphi(X_n)=\mE_\mu\varphi(X_0)$ for every $n\ge 0$. Due to the continuity of $\mathfrak{p}$ and the compactness of $\T$ the existence of a stationary distribution is ensured by standard techniques (see e.g. Corollary 4.18 on p. 31 in \cite{Hairer_06}) but its uniqueness is not so clear. We call a random walk $\textit{mixing}$ (or \textit{stable}) if there exists a measure $\mu_\ast$ such that
$$\mE_x \varphi(X_n) \to \int_\mathbb{T} \varphi d\mu_\ast \quad \textrm{as $n\to \infty$}$$
for every initial condition $x\in\mathbb{T}$ and $\varphi\in C(\mathbb{T})$, where $\mE_x$ is a conditional expectation provided the process is launched at $x\in\T$.

\subsection{Random walks in quasiperiodic environment}

Let $\mathfrak p : \mathbb{T} \rightarrow (0,1)$ be a continuous function, $\alpha$ an irrational angle. These along with a fixed point $x\in \mathbb T$ define a sequence $(p_j, q_j)$, $j\in\mathbb Z$, $p_j=\mathfrak{p}(x+j\alpha)$, $q_j=1-p_j$, that is called a \textit{quasiperiodic environment} defined by the triple $(\mathfrak p, \alpha, x)$. Each such environment gives rise to a nearest neighbour random walk $(\xi_n)$ on $\mathbb{Z}$ with $\xi_0=0$ and
\begin{equation}\label{E:1.3}
\mP_x( \xi_{n+1} = k+1 | \xi_n=k ) = p_k, \  \mP_x( \xi_{n+1} = k-1 | \xi_n=k ) = q_k, k\in \mathbb{Z}, n\ge 0.
\end{equation}
If $\mathfrak{p}$, $\alpha$ are fixed and $x\in \T$ is chosen randomly with uniform distribution, then $(p_j)_j$ becomes a stationary, ergodic process, which enables us to study (\ref{E:1.3}) using a well developed theory of random walks in random environment (see \cite{Zeitouni_04}). For example \cite{Zeitouni_04}, Theorem 2.1.2, implies that once $\mathfrak{p}$ and $\alpha$ are fixed and $\mathfrak{p}$ is \textit{symmetric}, i.e.
$$\int_{\T} \log \frac{\mathfrak{p}(x)}{\q(x)} dx = 0,$$
then $(\xi_n)$ is recurrent for Lebesgue almost every $x\in\T$. If the integral above is not zero then  $\mathfrak{p}$ is called \textit{asymmetric} and in that case $(\xi_n)$ is transient for Lebesgue almost every $x\in \T$.

One of important tools in the study of random walks in random environment is a certain auxiliary process called the environment viewed by the particle process. The general construction and some applications can be found in \cite{Zeitouni_04} (see e.g. the second proof of Theorem 2.1.9). This construction applied to (\ref{E:1.3}) leads to (\ref{E:1.1}), which points an interesting connection between these two processes.

\subsection{Motivation and results}

 A number $\alpha\not\in \mathbb{Q}$ is called \textit{Diophantine of type $(c, \tau)$}, $c>0$, $\tau\ge 0$, if
$$\bigg| \alpha - \frac{p}{q} \bigg| \ge \frac{c}{q^{2+\tau}} \quad \textrm{for every $p, q\in\mathbb{Z}, q\not = 0$.}$$
A number $\alpha$ is called \textit{Liouville} when it is not Diophantine of any type. It is a well known fact that for infinitely many natural $q$'s there exists $p$ such that
\begin{equation}\label{E:dioph}
\bigg| \alpha - \frac{p}{q} \bigg|<\frac{1}{q^2}.
\end{equation}
All such $q$'s are called \textit{close return times of $\alpha$}, while the numbers $d=| q\alpha - p|$ are called a \textit{close return distances of $\alpha$}.

Sinai \cite{Sinai_99} has proven the unique ergodicity, mixing and absolute continuity of stationary measure for (\ref{E:1.1}) and quenched local limit theorem for (\ref{E:1.3}) under the assumption that $\alpha$ is Diophantine and $\mathfrak{p}$ is sufficiently smooth. Kaloshin, Sinai \cite{Kaloshin_Sinai_00} have shown that the assumption that $\alpha$ is Diophantine can be dropped when $\mathfrak{p}$ is asymmetric. The unique ergodicity for symmetric $\mathfrak{p}$ and all irrational rotation numbers $\alpha$ has been obtained obtained by Conze and Guivarc'h \cite{Conze_Guivarch_00} under a mild assumption that $\mathfrak{p}$ is just continuous of bounded variation. An example of continuous $\mathfrak{p}$ (clearly, not of bounded variation) with infinitely many stationary distributions has been provided by J. Bremont in \cite{Bremont_99}.\footnote{There is a vast literature devoted to the study of limit properties of random walks in quasiperiodic environment (\ref{E:1.3}) and its relations to (\ref{E:1.1}). We are not going to discuss them in detail, an interested reader can check \cite{Bremont_01, Bremont_02, Bremont_09, Bremont_09B, DFS_21, Dolgopyat_Goldsheid_13, Dolgopyat_Goldsheid_18, Dolgopyat_Goldsheid_19, Dolgopyat_Goldsheid_20, Dolgopyat_Goldsheid_21, Goldsheid_07, Goldsheid_08}.}.

The first of our three results is the following.

\begin{Theorem}\label{T:1}
If $\alpha\not\in\mathbb{Q}$, $\mathfrak{p}$ is symmetric, $\log\frac{\mathfrak{p}(x)}{\q(x)}$ is continuous of bounded variation then the Markov process (\ref{E:1.1}) is uniquely ergodic.
\end{Theorem}

As mentioned before this theorem has been proven already in \cite{Conze_Guivarch_00}. The novelty is a new way of proving it stressing the relation to the random walk in quasiperiodic environment (\ref{E:1.3}). A crucial ingredient in the proof is so called strong ratio limit property (see \cite{Orey_61}).

Let $\mathcal{T}_k$ stands for the $k$-th ladder moment $\mathcal{T}_k = \min \{j \ge 0 : \xi_j = k \}$, $k\in \Z$. Theorem \ref{T:2} says the renewal theorem\footnote{For another relation to the renewal theorems see \cite{Kesten_77}, \cite{Lalley_86}.} for the ladder moments implies mixing of $X$.

\begin{Theorem}\label{T:2}
Let $\alpha\not\in\mathbb{Q}$, $\mathfrak{p}$ be symmetric, $\log\frac{\mathfrak{p}(x)}{\q(x)}$ continuous of bounded variation, $x\in \T$. If $\sum_{k\in \Z} \mP_x(\mathcal{T}_k=n) \to 0$ as $n\to \infty$, then $\mE_x\varphi(X_n) \to \int_{\T} \varphi d\mu_\ast$, where $\mu_\ast$ is the unique stationary measure.
\end{Theorem}

The renewal theorem above is unfortunately not easy to verify, and seems to be an interesting problem by itself. In the literature one can find some renewal theorems for independent non-identically distributed random variables (e.g. \cite{Kawata_56}) or for stationary ergodic processess (e.g. \cite{Berbee_79}, \cite{Lalley_86B}) but neither of these seems to apply. Nevertheless we are still able to show the following generic result.

\begin{Theorem}
\label{T:3}
Let $\alpha\not\in\mathbb{Q}$. Then there exists a dense $G_\delta$ subset in the space of zero mean absolutely continuous functions AC($\T$)\footnote{Zero mean refers to the Lebesgue measure. The space AC($\T$) is equipped with $\|\cdot\|_{BV}$ topology.} such that for each $f$ in this subset the Markov process (\ref{E:1.1}) with $\mathfrak{p}(x)= \frac{\exp f(x)}{1+\exp f(x)}$ is mixing.
\end{Theorem}

Theorem \ref{T:3} partially answers the question posed by D. Dolgopyat, B. Fayad and M. Saprykina in \cite{DFS_21} (Question 2.8 on p. 10). 

\subsection{Notation}

We are going to use symbols $\mP_x$ and $\mE_x$, $x\in \T$, referring to both the random walk (\ref{E:1.3}) on $\Z$ in the environment defined by $x\in \T$ and the process (\ref{E:1.1}) on $\T$ started at $x$. To avoid confusion, the process on $\Z$ is always denoted by $\xi$ and the process on the circle by $X$.

\subsection{Acknowledgements}

The work has been partially supported by Polish National Science Centre grant Preludium UMO-2019/35/N/ST1/02363 and partially by the European Union’s Horizon 2020 research and innovation programme under the Marie Skłodowska-Curie Grant Agreement No. 101034413.

\section{Basic facts about symmetric random walks on $\Z$}

Fix $\alpha$ and $\mathfrak{p}$ and consider a Markov chain (\ref{E:1.3}) on $\Z$ defined by the environment $x\in\T$. Let us define also a function $\mathfrak{m} : \Z \rightarrow \mathbb R$ by $\mathfrak{m}(0)=0$, $\mathfrak{m}(1)=1$, and
$$\mathfrak{m}(n)=1+\sum_{k=1}^{n-1}\prod_{j=1}^k \frac{\q(x+j\alpha)}{\mathfrak{p}(x+j\alpha)}=1+\sum_{k=0}^{n-2} \exp -S_kf(x+\alpha) \quad \textrm{for $n\ge 2$,}$$
and
$$\mathfrak{m}(-n)=-\sum_{k=0}^{n-1} \prod_{j=0}^k \frac{\mathfrak{p}(x-j\alpha)}{\q(x-j\alpha)}=-\sum_{k=0}^{-n+1} \exp S_kf(x) \quad \textrm{for $n\ge 1$,}$$
where $f(x)=\log\frac{\mathfrak{p}(x)}{\q(x)}$ and $S_kf$ denotes the $k$-th Birkhoff sum of $f$\footnote{If $k>0$ then $S_{-k}f(x)=f(x)+f(x-\alpha)+\cdots+f(x-k\alpha)$. Also $S_0f(x)=f(x)$.}. It is immediate to check that if $(\xi_n)$ evolves with the law (\ref{E:1.3}), then $\mathfrak{m}(\xi_n)$ is a martingale. As a consequence it can be proven that $(\xi_n)$ is recurrent if and only if $\mathfrak{m}(n)\to \pm\infty$ as $n\to \pm\infty$ (see e.g. (3.11) in \cite{DFS_21}, Section 3.2) and hence that $\xi$ with symmetric $\mathfrak{p}$ is recurrent not only for  Lebesgue almost every environment $x$ (as implied by already mentioned Theorem 2.1.2 \cite{Zeitouni_04}) but actually for every $x\in \T$.

\begin{Lemma}\label{L:recurrence}
If $\mathfrak{p}$ is symmetric, $f(x)=\log\frac{\mathfrak{p}(x)}{\q(x)}$ is continuous of bounded variation, then (\ref{E:1.3}) is recurrent for every $x\in \T$.
\end{Lemma}
\begin{proof}
It suffices to show $\mathfrak{m}(n)\to \pm\infty$ as $n\to \pm\infty$ for every $x$. Fix $n>0$ large and set $N=N(n)$ to be the number of positive close return times $q$ of $\alpha$ that are less than $n$. If $q$ is such close return time then $|S_q f(x)|<\var(f)$ by the Denjoy-Koksma inequality and the symmetry, thus $\mathfrak{m}(n)>N(n)\exp\var(f)$, which tends to infinity when $n\to \infty$. The proof when $n\to -\infty$ is an obvious modification.
\end{proof}

Let $\log\frac{\mfp(x)}{\mfq(x)}$ be continuous of bounded variation. Let $T_0=\min\{n > 0 : \xi_n = 0\}$ be the moment of the first return to $0$, and let $\rho_a(x)$, $a\in\Z\setminus\{0\}$, the expected number of visits in $a$ before $T_0$, i.e.
\begin{equation}\label{E:inv_measure_def}
\rho_{a}(x) = \mE_x\bigg( \sum_{j=0}^{T_0-1} \mathds{1}_{\{a\}}(\xi_j) \bigg), \quad a\in\mathbb{Z}.
\end{equation}
In our setting $\rho_a(x)$ can be computed explicitly in terms of $\mathfrak{p}$ and $\q$.
\begin{Lemma}\label{L:rho_computation}
Fix $x\in \T$. For a random walk (\ref{E:1.3}) in the environment $\mP_x$ we have
\begin{equation}\label{E:inv_measure_formula}
\rho_{a}(x)= \frac{\mathfrak{p}(x)}{\q(x+a\alpha)} \prod_{j=1}^{a} \frac{\mathfrak{p}(x+j\alpha)}{\q(x+j\alpha)}, \quad
\rho_{-a}(x)= \frac{\q(x)}{\mathfrak{p}(x-a\alpha)} \prod_{j=1}^{a} \frac{\q(x-j\alpha)}{\mathfrak{p}(x-j\alpha)},
\end{equation}
for $a>0$.
\end{Lemma}
\begin{proof}
Fix $x\in \T$. Let $a>0$ (for $a<0$ the proof is analogous), $T_a = \min\{n>0 : \xi_n=a\}$ be the moment of the first visit in $a$, $G=\{T_a<T_0\}$. The number of visits in $a$ before $T_0$ is clearly zero $\mP_x$-a.s. on $\{T_a\ge T_0\}$, therefore by the law of total probability and conditioning with respect to the stopping time $T_a$ we obtain
$$
\rho_a(x) = \mP_x(G) \mE_x\bigg( \sum_{j=0}^{T_0-1} \mathds{1}_{\{a\}}(\xi_j) \bigg| G \bigg) =  \mP_x(G) \mE_{x+a\alpha}\bigg( \sum_{j=0}^{T_0-1} \mathds{1}_{\{a\}}(\xi_j)\bigg).
$$
By (3.10), Section 3.2 in \cite{DFS_21} we have $\mP_x(G)=\mathfrak{p}(x)\frac{\mathfrak{m}(1)-\mathfrak{m}(0)}{\mathfrak{m}(a)-\mathfrak{m}(0)}$. Observe that $\sum_{j=0}^{T_0-1} \mathds{1}_{\{a\}}(\xi_j)$ has a geometric distribution with parameter equal to the $\mP_{x+a\alpha}$-probability that $T_0<T_a$, thus using again (3.10) in \cite{DFS_21} its expectation is
$$
\frac{\mathfrak{m}(a)-\mathfrak{m}(0)}{\q(x+a\alpha)\big(\mathfrak{m}(a)-\mathfrak{m}(a-1)\big)},
$$
Multiplying this by the value of $\mP_x(G)$ yields the assertion.
\end{proof}
It is worth pointing out an another easy-to-check relation of invariance of $\rho_a$ (see \cite{Dolgopyat_Goldsheid_19})
\begin{equation}\label{E:inv_measure}
\rho_{j+1}(x) \mathfrak{p}(x+j\alpha) + \rho_{j-1}(x) \q(x+j\alpha) = \rho_j(x) \quad x\in\T, j\in \Z.
\end{equation}
The recurrence and assumption that $\mathfrak{p}$, $\q$ are strictly positive implies the following statement called the strong ratio limit property (see \cite{Kingman_Orey_64} and \cite{Orey_61}, the proof can be found also in \cite{Freedman_83} Chapter 2.6 on p. 64).
\begin{Theorem}\label{T:SRLP}
Let $\mfp$ be symmetric such that $\log\frac{\mfp(x)}{\mfq(x)}$ is continuous of bounded variation. Then for every $\varepsilon>0$, $q\ge 2$ there exists $n_0$ such that
$$\bigg|\frac{\mP_x(\xi_n=a)}{\mP_x(\xi_n=b)} - \frac{\rho_{a}(x)}{\rho_{b}(x)}\bigg| < \varepsilon$$
for every $n\ge n_0$, $x\in\mathbb{T}$, $a,b \in [-q,q]$ both with the same parity as $n$.
\end{Theorem}
Since $(\xi_n)$ under $\mP_x$ is null recurrent the series $\mP_x(\xi_1=a)+\cdots+\mP_x(\xi_n=a)$ is divergent, the following corollary known as the Doeblin  ratio limit property \cite{Doeblin_38} is immediate.
\begin{Corollary}\label{C:SRLP}
Let $\mfp$ be symmetric such that $\log\frac{\mfp(x)}{\mfq(x)}$ is continuous of bounded variation. Then for every $\varepsilon>0$, $q\ge 2$ there exists $n_0$ such that
$$\bigg|\frac{\mP_x(\xi_1=a)+\cdots+\mP_x(\xi_n=a)}{\mP_x(\xi_1=b)+\cdots+\mP_x(\xi_n=b)} - \frac{\rho_{a}(x)}{\rho_{b}(x)}\bigg| < \varepsilon$$
for every $n\ge n_0$, $x\in\mathbb{T}$, $a,b \in [-q,q]$.
\end{Corollary}

\section{The proof of Theorem \ref{T:1}}
Throughout this section for any $z\in \T$, $j\in \Z$ we denote $z_j=z+j\alpha$. Let
$$
\nu_x^n = \frac{1}{M_n(x)} \sum_{j=0}^{n-1} \rho_j(x) \delta_{x_j}
$$
where $M_n(x)$ is the normalizing constant $M_n(x) = \sum_{j=0}^{n-1} \rho_j(x)$ and $\rho_j(x)$ were defined in (\ref{E:inv_measure_def}). Using (\ref{E:inv_measure_formula}) it is easy to verify the relation
\begin{equation}\label{E:rho_property}
\rho_{n+m}(x) = \rho_n(x)\rho_{m}(x_n), \quad \textrm{for $n, m \ge 0$, $x\in \T$.}
\end{equation}

\begin{Lemma}\label{L:1}
If $\log\frac{\mathfrak{p}(x)}{\q(x)}$ is continuous of bounded variation then
$$
\sup_{x\in\T}\max_{j=0,1,\ldots, n-1} \frac{\rho_j(x)}{M_n(x)} \to 0 \quad \textrm{as $n\to \infty$.}
$$
\end{Lemma}
\begin{proof}
 Let $f(x)=\log\frac{\mathfrak{p}(x)}{\q(x)}$, $C=\exp\var f$. For any $n>0$ let $N(n)$ be the number of close return times of $\alpha$ that are less than $n/2$. Take $j=0,1,\ldots, n-1$. If $j<n/2$ and $q$ is any close return time that is less than $n/2$ then $j+q$ is less than $n$, and thus $\rho_{j+q}(x)>C^{-1}\rho_j(x)$ by the Denjoy-Koksma inequality. Therefore $M_n(x)> N(n)C^{-1} \rho_j(x)$. The same inequality holds when $j\ge n/2$. The assertion follows since the Denjoy-Koksma inequality holds uniformly for all $x\in\T$ and the rate of divergence of $N(n)$ is independent of $x$.
\end{proof}

\begin{Proposition}\label{P:main}
There exists $L>0$ such that for every nonempty interval $I\subseteq \T$ there exists $q_0$ such that $\nu_x^q(I) \le L \nu_y^q(I)$ for every $x$, $y\in \T$ and every close return time $q$ of $\alpha$ with $q\ge q_0$.
\end{Proposition}
\begin{proof}
Since $\alpha$ is irrational, $\mu(I)>0$ for any stationary measure $\mu$ of the random walk $X$ (\ref{E:1.1}). By the Prokhorov theorem the set of all probability stationary measures is weakly-$\ast$ compact hence 
\begin{equation}\label{E:minimal_value}
\inf\mu(I)>0,
\end{equation}
where the infimum is taken over all probability stationary measures $\mu$. Denote $f(x)= \log\frac{\mathfrak{p}(x)}{\q(x)}$ and fix $x$, $y\in \T$. Let $q_m$ denote the $m$-th close return time, and let $d_m$ be the corresponding close return distance. Then $d_{m-1}\le 1/q_{m-1}$, $d_{m-1}>d_m$, and the points $x_j$, $j=0,1\ldots, q_m-1$, partition the circle into arcs of length at most $d_{m-1}+d_m<2d_{m-1}\le 2/q_{m-1}$, which shows that in each closed arc of length $2/q_{m-1}$ there is at least one point $x_j$, $j=0,1,\ldots q_m-1$, that belongs to this arc. Let $k$ be any index such that $x_k$ is not further to $y$ than $1/q_{m-1}$, and let $J$ be the arc joining $x_k$ and $y$. By standard arguments (e.g. Lemma 1.3 and the proof of Lemma 2.1 in \cite{deMelo_vanStrien_93}) the the subsets $J, R_\alpha(J), \ldots , R_\alpha^{q_m-k}(J)$ are pairwise disjoint and therefore
 $$
 \log\frac{\rho_j(x_k)}{\rho_j(y)} = \log\frac{\mathfrak{p}(x_k)}{\mathfrak{p}(y)}+\log\frac{\q(x_k+j\alpha)}{\q(y+j\alpha)} + \sum_{l=0}^{j-1} f(x_k+l\alpha) - f(y + l\alpha)
 $$
 for $j=1,\ldots, q_m-k$. Hence if $x_k$ and $y$ are sufficiently close (i.e. $q_m$ is sufficiently large) then
\begin{equation}\label{E:rho_property2}
\rho_j(x_k) \le C^2 \rho_{j}(y), \quad j=0,1,\ldots, q_m-k,
\end{equation}
where $C=\exp(\var(f))$. If $J'$ denotes the arc joining $x$ and $y_{k-q_m}$, then its length is at most $2/q_m$ and each point of the circle is in at most two subsets among  $J', R_\alpha(J'), \ldots , R_\alpha^{q_m-k}(J')$ (again Lemma 1.3 in \cite{deMelo_vanStrien_93}), hence by a similar argument as above
\begin{equation}\label{E:rho_property3}
\rho_j(x) \le C^3\rho_{j}(y_{q_m-k}), \quad j=0,1,\ldots, k.
\end{equation}
Moreover, by the Denjoy-Koksma inequality $C^{-1}<\rho_{q_m}(x)<C$, hence (\ref{E:rho_property}) and then (\ref{E:rho_property2}) give
\begin{equation}\label{E:rho_property4}
\rho_k(x)^{-1}\le C\rho_{q_m-k}(x_k)\le C^3 \rho_{q_m-k}(y).
\end{equation}
The relation (\ref{E:rho_property}) yields
$$
M_{q_m}(x) = M_k(x) + \rho_k(x)M_{q_m-k}(x_k), \quad M_{q_m}(y) = M_{q_m-k}(y) + \rho_{q_m-k}(y)M_k(y_{q_m-k}),
$$
therefore (\ref{E:rho_property2}), (\ref{E:rho_property3}) and (\ref{E:rho_property4}) imply
$$
M_{q_m}(x) \le C^3 \rho_k(x) \big(\rho_k(x)^{-1}M_{k}(y_{q_m-k}) + M_{q_m-k}(y) \big) 
$$
\begin{equation}\label{E:normalizing_property2}
\le C^4 \rho_k(x)\big( \rho_{q_m-k}(y)M_{k}(y_{q_m-k}) + M_{q_m-k}(y) \big)
=C^4\rho_k(x) M_{q_m}(y).
\end{equation}
Similarly
\begin{equation}\label{E:normalizing_property3}
M_{q_m}(y) \le C^4\rho_{q_m-k}(y) M_{q_m}(x).
\end{equation}
Combining (\ref{E:rho_property2}) and (\ref{E:normalizing_property2}) we have
\begin{equation}\label{E:2}
\frac{\rho_j(x)}{M_{q_m}(x)} \le C^6 \frac{\rho_{j-k}(y)}{M_{q_m}(y)} \quad \textrm{for $j=k,k+1,\ldots, q_m-1$,}
\end{equation}
and doing the same with  (\ref{E:rho_property3}) and (\ref{E:normalizing_property3}) we obtain
\begin{equation}\label{E:3}
\frac{\rho_j(x)}{M_{q_m}(x)} \le C^7 \frac{\rho_{q_m-k+j}(y)}{M_{q_m}(y)} \quad \textrm{for $j=0, 1,\ldots, k-1$.}
\end{equation}

Denote $q=q_m$ for short, fix a subinterval $I\subseteq \T$ and let $\Gamma_z = \{ j=0,1\ldots, q-1 : z_j\in I\}$, $z\in\T$. Then by (\ref{E:2}) and (\ref{E:3}) we get
$$
\nu_x^q (I) = \sum_{
\substack{
j\in \Gamma_x \\
j<k
}} \frac{\rho_j(x)}{M_q(x)}
+\sum_{
\substack{
j\in \Gamma_x \\
j\ge k
}} \frac{\rho_j(x)}{M_q(x)}
\le 
C^7\sum_{
\substack{
j\in \Gamma_x \\
j<k
}}\frac{\rho_{q-k+j}(y)}{M_q(y)}
+C^6\sum_{
\substack{
j\in \Gamma_x \\
j\ge k
}} \frac{\rho_{j-k}(y)}{M_q(y)}.
$$
Since $J, R_\alpha(J), \ldots, R^{q-k}(J)$ are pairwise disjoint there are at most two $j$'s with $j<k$ such that $j\in \Gamma_x$ and $q-k+j\not\in \Gamma_y$ or $j\not\in \Gamma_x$ and $q-k+j\in \Gamma_y$. For the same reason there are at most four $j$'s, $j\ge k$ such that $j\in \Gamma_x$ and $j-k\not\in \Gamma_y$ or $j\not\in \Gamma_x$ and $j-k\in \Gamma_y$. Therefore we can write
$$
\nu_x^q(I) \le C^7\sum_{
\substack{
j\in \Gamma_y \\
j<k
}}\frac{\rho_{q-k+j}(y)}{M_q(y)}
+ C^6\sum_{
\substack{
j\in \Gamma_y \\
j\ge k
}} \frac{\rho_{j-k}(y)}{M_q(y)} + 4C^7\sup_{z\in \T}\max_{j=0,1\ldots, q-1} \frac{\rho_j(z)}{M_q(z)}.
$$
The sum of the first two terms is less than $C^7 \nu_y^q(I)$. The supremum tends to zero by Lemma \ref{L:1}, thus it is smaller than  $\frac 1 4 \nu_y^q(I)$ for $q$ sufficiently large by  (\ref{E:minimal_value}). This completes the proof with $L=2C^7$.
\end{proof}

\begin{proof}[Proof of Theorem \ref{T:1}]
We are going to show there exists a stationary measure $\mu_\ast$ with the property that every ergodic stationary measure is absolutely continuous with respect to $\mu_\ast$. This immediately implies the assertion.

Let $\mu_\ast$ be an arbitrary accumulation point of an arbitrary sequence of the form $\nu_{y}^{q_n}$, where $y\in \T$, $q_n\to\infty$ are close return times of $\alpha$. Let $\mu$ be an arbitrary ergodic stationary measure, let $x\in \T$ be a generic point of this measure, and let $\varphi$ be the characteristic function of some arc $I\subseteq \T$. The goal is to show that $\mu(I)\le L\mu_\ast(I)$, where $L$ is the constant from Proposition \ref{P:main}. 

Let $q=q_n$ for some $n$, $\mathcal{I}_s=[sq,(s+1)q)$, $s\in\Z$, and let $\varphi_s(j) = \varphi(x+j\alpha)$ for $j\in \mathcal{I}_s$ and $0$ otherwise. Let $\tau_s$ be the moment of the first visit in $\mathcal{I}_s$, $s\in\Z$. Then
$$
\mE_x(\varphi(X_1) + \cdots + \varphi(X_n))
= \sum_{s\in\Z} \mE_x  \mE_{\xi_{\tau_s}}\big( \varphi_s(\xi_{0}) + \cdots + \varphi_s(\xi_{n-\tau_s})\big).
$$

Let $\varepsilon>0$, $n_0$ given in Corollary \ref{C:SRLP}. Observe that $\sum_{s\in\Z} \mE_x\mE_{\xi_{\tau_s}}\big[ \mathds{1}_{\mathcal{I}_s}(\xi_{0})+\cdots \mathds{1}_{\mathcal{I}_s}(\xi_{n-\tau_s}) \big]=n$ and $\sum_{s\in\Z} \mP_x(\tau_s\ge n-n_0) \le \sum_{s\in\Z} \sum_{j=0}^{n_0} \mP_x(\xi_{n-j} \in \mathcal{I}_s) = n_0$, thus
$$\sum_{s\in\Z} \mE_x \mathds{1}_{\{\tau_s\ge n-n_0\}} \mE_{\xi_{\tau_s}}\big[ \mathds{1}_{\mathcal{I}_s}(\xi_{0})+\cdots \mathds{1}_{\mathcal{I}_s}(\xi_{n-\tau_s}) \big] \le n_0^2
$$
and
\begin{equation}\label{E:limit1}
\frac 1 n \sum_{s\in\Z} \mE_x \mathds{1}_{\{\tau_s < n-n_0\}} \mE_{\xi_{\tau_s}}\big[ \mathds{1}_{\mathcal{I}_s}(\xi_{0})+\cdots \mathds{1}_{\mathcal{I}_s}(\xi_{n-\tau_s}) \big] \to 1 \quad \textrm{as $n\to \infty$}.
\end{equation}
By Corollary \ref{C:SRLP}
$$
 \mE_{\xi_{\tau_s}}\big( \varphi_s(\xi_{0}) + \cdots + \varphi_s(\xi_{n-\tau_s})\big) 
 \le \big( \nu_{x}^q(I) + \varepsilon \big) \mE_{\xi_{\tau_s}}\big[ \mathds{1}_{\mathcal{I}_s}(\xi_{0})+\cdots \mathds{1}_{\mathcal{I}_s}(\xi_{n-\tau_s}) \big]
$$
$\mP_x$ a.s. on $\{\tau_s< n-n_0\}$. By this, (\ref{E:limit1}) and Proposition \ref{P:main} for $n$ large we have
$$
\frac{\mE_x(\varphi(X_1) + \cdots + \varphi(X_n))}{n}
=
\sum_{s\in\Z} \frac{\mE_x  \mE_{\xi_{\tau_s}}\big( \varphi_s(\xi_{0}) + \cdots + \varphi_s(\xi_{n-\tau_s})\big)}{n} \le L \nu_{y}^q(I) + 2\varepsilon.
$$
The left hand side tends to $\mu(I)$ by the Birkhoff ergodic theorem, thus taking $\varepsilon\to 0$ and then $q\to \infty$ completes the proof.
\end{proof}

\begin{Remark}\label{Remark}
For a moment let us fix three functions $\mathfrak{p}$, $\q$, $\mathfrak{s} : \T \to \mathbb{R}_+$ with a $\mathfrak{p}(x)+\q(x)+\mathfrak{s}(x)=1$ for every $x\in \T$, and consider a Markov process on $\T$ with the transition kernel
\begin{equation}\label{E:walk_def_new}
p(x, \cdot) = \mathfrak{p}(x) \delta_{x+\alpha} + \mathfrak{s}(x)\delta_x + \q(x) \delta_{x-\alpha}, \quad x\in \T.
\end{equation}
It gives rise to a random walk $(\xi_n)$ on $\Z$ in the same way as (\ref{E:1.3}). The functions $M$ and $\rho$ defined in Section 2 preserve their meaning with exactly the same formulae, and therefore the proof of Theorem \ref{T:1} works for (\ref{E:walk_def_new}) without any changes if only $\log\frac{\mathfrak{p}(x)}{\q(x)}$ is continuous of bounded variation.

Let us now go back to the original process (\ref{E:1.1}). Then $(X_{2n})$ is of the form (\ref{E:walk_def_new}), and $\log\frac{\mathfrak{p}(x)\mathfrak{p}(x+\alpha)}{\q(x)\q(x-\alpha)}$ is of bounded variation if $\log\frac{\mathfrak{p}(x)}{\q(x)}$ is of bounded variation. Therefore not only $(X_n)$ is uniquely ergodic but also $(X_{2n})$ is uniquely ergodic under the same assumptions. This is not a consequence of Theorem \ref{T:1}.
\end{Remark}

\section{The proof of Theorem \ref{T:2}}

 It is more convenient to show the assertion for $(X_{2n})$ and then use the conditional expectation to obtain the result for $(X_n)$.
 
Fix $\varphi\in C(\T)$ and $\varepsilon>0$. Let $(x_n)$ be an arbitrary sequence of points in $\T$, and let $(q_n)\subseteq \mathbb{N}$ be an arbitrary sequence of close return times of $\alpha$ with $q_n\to \infty$. Lemma \ref{L:1} and the relation (\ref{E:inv_measure}) applied twice imply that every accummulation point of the sequence of measures
 $$\nu_{x_n}^{q_n} = \frac{1}{M_{q_n}(x_n)} \sum_{j=0}^{q_n} \rho_{2j}(x_n)\delta_{x_n+2j\alpha},$$
where $M_{q_n}(x)$ is the normalizing constant, is some stationary measure of $(X_{2n})$. By Remark \ref{Remark} $(X_{2n})$ is uniquely ergodic with the same stationary measure $\mu$ as $(X_n)$, thus $\nu_{x_n}^{q_n}$ converge to $\mu$. Moreover, a similar argument as in Proposition 4.1.13 \cite{Katok_Hasselblatt_95} shows the convergence is actually uniform over $x$, thus there exists $q$ so large that for any $x\in\T$
\begin{equation}\label{E:measure_close}
\bigg| \int_{\T} \varphi d\nu_x^q - \int_{\T} \varphi d\mu \bigg| < \varepsilon.
\end{equation}
Let $\mathcal{I}_s=[qs, q(s+1))$, $s\in\mathbb{Z}$. The ladder moment $\mathcal{T}_{sq}$ is the moment of the first visit in $\mathcal{I}_s$ if $s>0$, while $\mathcal{T}_{(s+1)q-1}$ is the moment of the first visit in $\mathcal{I}_s$ for $s<0$. Fix $x\in\T$. Theorem \ref{T:SRLP} and (\ref{E:measure_close}) give $n_0$ so large that for $s>0$ we have

\begin{equation}\label{E:4.1}
\bigg| \mE_{\xi_{\mathcal{T}_{sq}}}  \varphi( x+\xi_{n-\mathcal{T}_{sq}} \alpha )   - \int_{\T} \varphi d\mu \bigg| < 2\varepsilon \quad \textrm{a.s. on $\{ \mathcal{T}_{sq} \le n-n_0 \}$.}
\end{equation}
Since the renewal theorem for the ladder moments holds by assumption we get
\begin{equation}\label{E:maximum}
\sum_{s\in\Z} \mP_x \big(\mathcal{T}_{sq} \in (n-n_0, n]    \big)\le \sum_{j=0}^{n_0-1}\sum_{s\in\Z} \mP_x \big( \mathcal{T}_s = n-j \big) \to 0 \quad \textrm{as $n\to \infty$.}
\end{equation}
We have now
$$\mE_x \varphi( x+\xi_n \alpha ) = \sum_{s\in\Z} \mE_x \mathds{1}_{\{ \mathcal{T}_{sq} > n-n_0 \}} \mE_{\xi_{\mathcal{T}_{sq}}}  \varphi( x+\xi_{n-\mathcal{T}_{sq}} \alpha ) $$
$$
+ \sum_{s\in\Z} \mE_x \mathds{1}_{\{ \mathcal{T}_{sq} \le n-n_0 \}} \mE_{\xi_{\mathcal{T}_{sq}}}  \varphi( x+\xi_{n-\mathcal{T}_{sq}} \alpha ) $$
By (\ref{E:maximum}) the first sum is bounded by $\varepsilon$ if $n$ is sufficiently large. By (\ref{E:4.1}) the second sum is closer than $2\varepsilon$ to $\int_\T \varphi d\mu$.

\section{Proof of Theorem \ref{T:3}}

First let us make an easy observation that the renewal theorem for the ladder moments $\mathcal{T}_s$ with respect to $\mP_x$, $x\in \T$, is implied if
$$
\mP_x \bigg( \max_{j\le n} \xi_j = \xi_n \bigg) \to 0
\quad \textrm{and} \quad
\mP_x \bigg( \min_{j\le n} \xi_j = \xi_n \bigg) \to 0 \quad \textrm{when $n\to \infty$}.
$$

Since the intersection of two dense $G_\delta$ sets is a dense $G_\delta$ it suffices to show that $\mP_x \bigg( \max_{j\le n} \xi_j = \xi_n \bigg) \to 0$ as $n\to \infty$ for a generic $\mathfrak{p}$ (for the minimum it is obtained by the symmetry). Let $\mathcal{P}$ denote the space of all absolutely continuous functions with zero mean with respect to the Lebesgue measure. Put
$$\mathcal{A}_m = \bigg\{ \mathfrak{p} \in \mathcal{P} : \exists_{n\ge 1} \forall_{x\in \T} \mP_x\bigg(\max_{j\le n} \xi_{j} = \xi_n \bigg) < \frac{1}{m} \bigg\}.$$
Clearly $\mathcal{A}_m$ is open for every $m$.

\begin{Lemma}
If $\mathfrak{p} \in \bigcap_{m=1}^\infty \mathcal{A}_m$, then $\mP_x (\max_{j\le n} \xi_{j} = \xi_n ) \to 0$ for every $x\in\T$ as $n\to\infty$.
\end{Lemma}
\begin{proof}[Proof of the Lemma]
Let $m_0$ be natural. Since $\mathfrak{p}\in\mathcal{A}_{m_0}$, there is $n_0$ be such that $\mP_x \bigg(\max_{j\le {n_0}} \xi_{j} = \xi_{n_0} \bigg) < \frac{1}{m_0}$ for each $x\in\T$. Then, for fixed $x\in \T$ and $n\ge n_0$,
$$\mP_x \bigg(\max_{j\le n} \xi_{j} = \xi_n \bigg) \le \mP_x \bigg(\max_{j\le n_0} \xi_{n-j} = \xi_n \bigg) = \mE_x \mP_{\xi_{n-n_0}} \bigg(  \max_{j\le n_0} \xi_{j} = \xi_{n_0} \bigg) < \frac{1}{m_0}.$$
\end{proof}

 Trigonometric polynomials form a dense subset of $\mathcal P$ equipped with $\|\cdot \|_{BV}$ topology. Thus what remains is to show is that for every zero mean trigonometric polynomial $f$ the function $\mathfrak{p}(x)=\frac{\exp f(x)}{1+\exp f(x)}$ is an element of each  $\mathcal{A}_m$, $m\ge 1$.
 
 Take $\delta>0$ small and $n$ large. Clearly
$$\mP_x \bigg(\max_{j\le n} \xi_j - \xi_{n} = 0 \bigg) \le \mP_x \bigg(\max_{j\le n} \frac{\xi_{j}}{\sqrt{n}} -\frac{ \xi_{n}}{\sqrt{n}} \le \delta \bigg)$$
for every $x\in \T$. Since $f$ is a trigonometric polynomial, the assumptions of Theorem 6.1 and Corollary 6.2 in \cite{Dolgopyat_Goldsheid_21} are satisfied and the process $t\longmapsto \frac{\xi_{\lfloor nt \rfloor}}{\sqrt{n}}$, $t\in [0,1]$, converges in distribution to the Brownian motion $B_t$ defined on some space $(\Omega, \mathcal{F}, \mP)$ with variance $\sigma^2(x)t>0$, $t\in [0,1]$, that is uniformly bounded for $x\in \T$. Thus
$$\mP_x \bigg(\max_{j\le n} \frac{\xi_{j}}{\sqrt{n}} -\frac{ \xi_{n}}{\sqrt{n}} \le  \delta \bigg) \le \mP_x \bigg(\max_{j\le k} \frac{\xi_{j \lfloor \frac{n}{k} \rfloor }}{\sqrt{n}} -\frac{ \xi_{n}}{\sqrt{n}} \le \delta \bigg) \to  \mP \bigg(\max_{j \le k} B_{j/k} - B_1 \le \delta \bigg)
$$
uniformly over $x\in T$. But Proposition 11.13 on p. 206 in \cite{Kallenberg_97} implies that
$$\mP\bigg( \max_{s\le 1} B_s - B_1 \le \delta \bigg) = 2\int_0^\delta \frac{1}{\sigma(x)\sqrt{2\mathfrak{p}i} } e^{-u^2/2\sigma(x)^2}du.$$
The observation that
$$
\mP \bigg( \max_{s\le 1} B_s - \max_{j\le k} B_{j/k} > \delta \bigg) \le \delta
$$
if $k$ is sufficiently large, leads to
$$
    \mP \bigg(\max_{j \le k} B_{j/k} - B_1 \le \delta \bigg)<2\int_0^\delta \frac{1}{\sigma(x)\sqrt{2\mathfrak{p}i} } e^{-u^2/2\sigma(x)^2}du + 2\delta, $$
which can be made arbitrarily small. This completes the proof.

\appendix
\section{}

Theorem \ref{T:SRLP} is slightly stronger than \cite{Orey_61} since it contains additional information about the uniform convergence over $x\in\T$. Here we present how to modify the proof presented in \cite{Freedman_83} Chapter 2.6 on p. 64 to get the necessary statement.

Assume contrary to the claim that there exist $\varepsilon>0$, $m_0\in\Z$ and sequences $(x_n)$, $(k_n)$ with
$$
\bigg| \frac{\mP_{x_n}\big(\xi_{k_n}=m_0\big)}{\mP_{x_n}\big(\xi_{k_n} = 0\big)} - \rho_{x_n}(m_0) \bigg| > \varepsilon.
$$
By the compactness of $\T$ we can assume that $x_n \to x$. By diagonal argument we can find a subsequence of $(k_n)$ (denoted by the same symbol) with
$$
\frac{\mP_{x_n}\big(\xi_{k_n} = m \big)}{\mP_{x_n}\big(\xi_{k_n} = 0\big)} \to \widetilde{\rho}_x(m)
$$
for all $m\in \Z$, where $\widetilde{\rho}_x(m_0) \not= \rho_x(m_0)$. We have
$$
\tilderho(m-1)\mathfrak{p}(x+(k-1)\alpha) + \tilderho(m+1)\q(x+(k+1)\alpha)
$$
$$
= \lim_{n \to \infty} \frac{\mP_{x_n}\big(\xi_{k_n} = m-1 \big)}{\mP_{x_n}\big(\xi_{k_n} = 0\big)} \mathfrak{p}(x+(k-1)\alpha)  + \frac{\mP_{x_n}\big(\xi_{k_n} = m+1 \big)}{\mP_{x_n}\big(\xi_{k_n} = 0\big)}\q(x+(k+1)\alpha)
$$
$$
= \lim_{n \to \infty} \frac{\mP_{x_n}\big(\xi_{k_n+1} = m\big)}{\mP_{x_n}\big(\xi_{k_n} = 0\big)}.
$$
The above converges to $\tilderho(m)$ if only the assertion of Lemma 55 in \cite{Freedman_83} holds uniformly over $x\in\T$. This yields the relation 
$$
\tilderho(m-1)\mathfrak{p}(x+(k-1)\alpha) + \tilderho(m+1)\q(x+(k+1)\alpha) =\tilderho(m),
$$
which is the invariant measure relation (\ref{E:inv_measure}) that uniquely determines $\rho_x$ (see \cite{Derman_54}). This contradicts $\tilderho(m_0)\not= \rho(m_0)$.

We are left with showing that Lemma 55 \cite{Freedman_83} holds uniformly over $x\in\T$. This is done by a simple observation that $\varepsilon$ in Lemma 55 can be chosen uniformly over $x\in\T$ by the compactness of $\T$ and the following Lemma that should replace (59).

\begin{Lemma}
If $\log\frac{\mathfrak{p}(x)}{\q(x)}$ is continuous of bounded variation, $\mathfrak{p}$ symmetric, then for every $t \in (0,1)$, $\varepsilon>0$ there exists $n_0$ such that
$$ \frac{t^n}{\mP_x(\zeta_n=0)} < \varepsilon$$
for every $x\in\T$ and $n\ge n_0$.
\end{Lemma}
\begin{proof}
It is sufficient to show that $\frac{-\log \mP_x (\xi_n = 0)}{n}$ converges to $0$ uniformly over $x\in \T$. Clearly the limit inferior of this sequence is $0$ pointwise since otherwise the Borel-Cantelli lemma contradicts the recurrence. The sequence $-\log \mP_x (\xi_n = 0)$ is subadditive, and what remains is to check that the proof of Fekete's lemma (Lemma 49 \cite{Freedman_83}) gives uniform convergence.

Take $\varepsilon>0$. By the compactness of $\T$ there exist natural numbers $N_1,\ldots, N_r$ and open neighbourhoods $V_{1}, \cdots, V_{r}$ covering the circle such that
$$
\frac{-\log \mP_z (\xi_{N_j} = 0)}{N_j} < -\log(\varepsilon) \quad \textrm{for $z \in V_j$, $j=1,\ldots r$.}
$$
By subadditivity, for any $j$, $k\ge 1$, $s=0,1,\ldots, N_j-1$ and $z\in V_j$ we have
$$
\frac{-\log\mP_z(\xi_{kN_j+s}=0)}{kN_j+s}\le \frac{kN_j}{kN_j+s}\cdot \frac{-\log\mP_z(\xi_{kN_j}=0)_{kN_j}-\log\mP_z(\xi_{s}=0)}{kN_j}
$$
$$
\le -\log(\varepsilon) + \frac{-\log\mP_z(\xi_{s}=0)}{kN_j}.
$$
Since $s$ is bounded by the largest number among $N_j$'s, the second term above is less than $-\log\varepsilon$ if only $k$ is sufficiently large. This completes the proof.
\end{proof}

\bibliographystyle{plain}
\bibliography{Mixing_generic}

\end{document}